\documentclass{article} 
 
\usepackage{ fullpage}  
\usepackage{amsmath, amssymb, amsfonts, amsthm, bbm, mathtools}  
 
\usepackage[normalem]{ulem} 
\usepackage{color,enumitem} 

\usepackage{thmtools}
\usepackage{thm-restate}

\newcommand{\red}[1]{{\color{red} #1}}

\newcommand{\comment}[1]{}

\newtheorem*{theorem*}{Theorem}

\newtheorem{lemma}{Lemma}
\newtheorem{corollary}{Corollary}

\newtheorem{definition}{Definition}

\newtheorem{example}{Example}
\newtheorem*{example*}{Example}

\newcommand{\bvec}[1]{\pmb{#1}}

\newcommand{\Real}{\mathbb{R}}
\newcommand{\Natural}{\mathbb{N}}

\newcommand{\good}{valid}

\def\avoid{\mathrm{avoid}}
\def\Ext{\mathrm{Ext}}
\def\ter{\mathrm{t}}
\newcommand{\RR}{L}
\newcommand{\zz}{\widetilde{w}}

\newcommand{\ratio}{\mathrm{ratio}}

\date{}

\title{The Lov\'{a}sz Local Lemma is Not About Probability}

\author{
Dimitris Achlioptas\\
University of Athens
\and
Kostas Zampetakis\\
UC 
Santa Cruz
}

\begin{document}

\maketitle

\begin{abstract}
Given a collection of independent events each of which has strictly positive probability, the probability that all of them occur is also strictly positive. The Lov\'asz local lemma (LLL), a mainstay of probabilistic combinatorics for more than 40 years, asserts that this remains true if the events are not too strongly negatively correlated. The formulation of the lemma involves a graph with one vertex per event, with edges indicating potential negative dependence. In its original formulation and its many subsequent extensions and refinements, the ``Local" in LLL reflects that the condition for the negative correlation to not be too great can be expressed solely in terms of the neighborhood of each vertex. In contrast to this local view, Shearer developed an exact criterion for the avoidance probability to be strictly positive, but it involves summing over all independent sets of the graph. 

In this work we make two contributions.  The first is to develop a hierarchy of increasingly powerful, increasingly non-local lemmata for bounding the avoidance probability from below, each lemma associated with a different set of walks in the graph. The LLL is at the bottom of this hierarchy, corresponding to the set of all possible walks on the graph, while at the top of the hierarchy is Shearer's exact criterion, corresponding to a set of self-avoiding walks on the graph intimately related with the computation tree of the independent set polynomial. Already, at its second level our hierarchy is stronger than all known local lemmata. To demonstrate the power of the hierarchy we use it to prove new bounds for the negative-fugacity singularity of the hard-core model on several lattices, a central problem in statistical physics. Our bounds match the conjectured values to three decimal digits with minimal computational effort.

Our second contribution is to prove that Shearer's connection between the probabilistic setting and the independent set polynomial holds for \emph{arbitrary supermodular} functions, not just probability measures. 
This means that all LLL machinery developed so far, including our new hierarchy, can be employed to bound from below an arbitrary supermodular function, based only on information regarding its value at singleton sets and partial information regarding their interactions. To demonstrate the power of this, we show that it trivially implies both the quantum LLL of Ambainis, Kempe, and Sattath~[JACM 2012], and the  quantum Shearer criterion of Sattath, Morampudi, Laumann, and Moessner~[PNAS 2016], where in both results the supermodular function bounded is relative subspace dimension.
\end{abstract}

\newpage

\section{Introduction}\label{subsec:PMLLL}

In the probabilistic method, introduced by Erd\H{o}s in his seminal 1947 paper~\cite{PM}, one proves that a certain object exists by demonstrating a probability distribution under which it has strictly positive probability. Much of the power of the method stems from the fact that any, arbitrarily bad, multiplicative approximation of the object's probability suffices: if it was strictly positive to begin with, it will remain strictly positive even after reduction by an arbitrarily large factor. For example, if we can prove that the object's probability is at least  the probability that a number of \emph{independent} trials all have a favorable outcome, we win.

The Lov\'{a}sz Local Lemma (LLL), a mainstay of the probabilistic method and a workhorse of combinatorics and theoretical computer science, works in the manner suggested above: it bounds from below the probability that a complex experiment has a favorable outcome, by the probability that a number of \emph{independent} trials all have a favorable outcome. More specifically, given a universe of candidate objects $\Omega$, the (potentially empty) set of desired objects is defined as the complement of $F_1 \cup \cdots \cup F_n$, each set $F_i$ containing undesirable objects that share some ``flaw." Thus, to establish the existence of the desired (flawless) objects, it suffices to establish that $\mu(\overline{F}_1 \cap \cdots \cap \overline{F_n})>0$ for some probability measure $\mu$ on $\Omega$ of our choice. The LLL enables this by offering a condition under which $\mu(\overline{F}_1 \cap \cdots \cap \overline{F_n}) \ge \prod_{i \in [n]} (1-r_i)>0$, where $\mu(F_i) \le r_i < 1$.

To express the LLL condition, write $[n] = \{1,2,\ldots,n\}$ and $p_i = \mu(F_i)$. Clearly, if the sets $\{F_i\}_{i\in [n]}$ are allowed to intersect arbitrarily, the best possible lower bound for $\mu(\overline{F}_1 \cap \cdots \cap \overline{F_n})$ is $1 - \sum_i p_i$, as the sets, in principle, could be disjoint. To improve upon the union bound we need to constrain the flaw overlaps.  A natural way to express this restriction is via a graph $G$ on $[n]$. Specifically, letting $\Gamma_i=\Gamma_i(G)$ denote the neighborhood of vertex $i$ in $G$, we say that $G$ is a \emph{dependency graph} for $\{F_i\}_{i=1}^n$ with respect to $\mu$, if
 for every $i \in [n]$ and every $\{j_1, j_2, \ldots \} \subseteq [n] \setminus (\Gamma_i \cup \{i\})$, conditioning on the avoidance of non-neighboring flaws does not change the probability of $F_i$, i.e., $\mu(F_i \mid \overline{F_{j_1}} \cap \overline{F_{j_2}} \cap \cdots) = p_i$. So, at the extremes, a complete dependency graph (clique) conveys no information at all about the conditional probabilities, while an empty dependency graph conveys that the $n$ events are mutually independent. 
Given $\bvec{p}, G$, we let $\avoid(\bvec{p}, G)$ denote the smallest possible value of the \emph{avoidance probability} $\mu(\overline{F}_1 \cap \cdots \cap \overline{F_n})$, when $\mu(F_i) = p_i$ and $G$ is a dependency graph.

\paragraph{The Lov\'asz Local Lemma.}
The original, symmetric formulation~\cite{erdHos1973problems} of the Lov\'{a}sz Local Lemma asserts that if $p_i \le p$, $|\Gamma_i| \le d$, and $pd \le 1/4$, then $\avoid(\bvec{p}, G)>0$. The general, asymmetric LLL formulation~\cite{erdHos1991lopsided} asserts that if there exist $r_1, \ldots, r_n \in [0,1)$ such that $p_i \le r_i \prod_{j \in \Gamma_i} (1-r_j)$, for all $i \in [n]$, then $\avoid(\bvec{p}, G) \ge \prod_{i\in[n]} (1-r_i)$. 

\paragraph{Shearer's Criterion.}
While the LLL gives an easily testable sufficient condition for $\avoid(\bvec{p}, G)>0$, a seminal result of Shearer~\cite{shearer} gave an \emph{exact} characterization, connecting $\avoid(\bvec{p}, G)$ to the independent set polynomial of the dependency graph $G$.  Recall that the independent set polynomial of graph $G$ on variables $x_1, \ldots, x_n$ is  $Z_G(\bvec{x}) = \sum_{I \in \mathrm{Ind}(G)} \prod_{i\in I} x_i$, where $\mathrm{Ind}(G)$ is the set of all independent sets of $G$. 
Shearer~\cite{shearer} proved the first equivalence below, and Scott and Sokal~\cite{SS} the second.
\[
\avoid(\bvec{p}, G) > 0 \quad \Longleftrightarrow \quad Z_G(-\bvec{p} ; S) > 0 \text{ for every $S \subseteq [n]$ } \quad \Longleftrightarrow\quad  Z_G(- \lambda \bvec{p}) > 0 \text{ for every $\lambda \in [0,1]$} \enspace .
\]
Let $\mathcal{S}(G) =\{\bvec{p} \in [0,1)^n: \avoid(\bvec{p}, G) > 0\}$. Shearer, also proved that if $\bvec{p} \in \mathcal{S}(G)$ then $\avoid(\bvec{p}, G) = Z_G(-\bvec{p})$. Naturally, verifying Shearer's condition is, in general, intractable as it involves summing over all independent sets in $G$. In light of Shearer's criterion, one can think of the LLL as asserting the following, which was independently established by Dobrushin~\cite{Dobrushin1996EstimatesOS, Dobrushin1996} in the context of statistical physics, 20 years after the LLL.

\paragraph{Dobrushin's Condition.}
If there exist $r_1, \ldots, r_n \in [0,1)$ such that $p_i \le r_i \prod_{j \in \Gamma_i} (1-r_j)$, for all $i \in [n]$, then $Z_G(-\bvec{p}; S) \ge \prod_{i\in S} (1-r_i) > 0$, for every $S \subseteq [n]$. 

\subsection{Our Contribution}

We make two contributions. \smallskip

Our first contribution is a hierarchy of sufficient conditions for the avoidance probability to be strictly positive which interpolates smoothly between the LLL and Shearer's exact criterion. Our hierarchy affords a tremendous amount of flexibility to the user, in order to adapt the condition to the problem at hand. The main innovation is that it has a strongly \emph{graphical} flavor, going (arbitrarily) far beyond vertex degree considerations. At the same time, as soon as one takes into account the subgraph induced by each inclusive neighborhood $\Gamma_i \cup \{i\}$ it already dominates every known \emph{local} condition, e.g., the asymmetric LLL~\cite{erdHos1991lopsided}, the cluster-expansion LLL~\cite{cluster}, and the non-backtracking LLL~\cite{ICALPIS}. 

To demonstrate the power of our hierarchy, we give new lower bounds for the negative-fugacity singularity of the hard-core model on several lattices, a central problem in statistical physics \cite{Gaunt1965HardSphereLG, todo1999transfer, chan2012series}. Our bound of $-0.1191$ for $\mathbb{Z}^2$ improves upon the previous best rigorous lower bound of $-0.113$ by Kolipaka, Szegedy, and Xu~\cite{cliqueLLL} and matches the conjectured~\cite{todo1999transfer} value $-0.11933888188...$, to three decimal digits. We achieve a similar level of accuracy relative to the conjectured values for other lattices. As we will see, each bound is derived by selecting an integer $q$ and enumerating all walks starting at the origin of the lattice that have length at most $q$ and satisfy some additional conditions. The results presented here come from relatively small values of $q$, corresponding to computational effort in the order of minutes on a laptop computer. We plan to report improved bounds in subsequent versions of this work and we \emph{conjecture} that for every $\epsilon >0$ there is $q=q(\epsilon)$ which brings our lower bound within $\epsilon$ of the exact singularity.\medskip

Our second contribution is to show that both the Lov\'{a}sz Local Lemma and Shearer's criterion hold for \emph{arbitrary supermodular} functions, not just probability measures. Specifically, we show how to: 
\begin{itemize}
\item
Lower bound an arbitrary supermodular function over its entire domain, in terms of (a) its value at singleton sets and (b) a graph $G$ that expresses the interactions between these sets, and,
\item
Express the optimal lower bound implied by this information in terms of the independent set polynomial of the graph $G$.
\end{itemize} 
As an example of the power afforded by the elevation to supermodular functions, we show that it readily implies both the original quantum LLL of Ambainis, Kempe, and Sattath~\cite{QLLL}, and the more recent quantum Shearer-like criterion of Sattath, Morampudi, Laumann, and Moessner~\cite{SPNAS}. In both cases, the function bounded from below is the relative dimension.

\section{Statement of Results}
Throughout the paper we will be referring to a graph $G$ with vertex set $[n]$ and edge set $E$. We denote the subgraph of $G$ induced by a set of vertices $S \subseteq [n]$ by $G[S]$.

\subsection{Shearer's Condition in Terms of Self-bounding Walks}

A \emph{walk} $w$ on $G$ is a string in $[n]^*$ where every two successive characters are adjacent in $G$. 
\begin{definition}[Self-bounding walks]
A walk on a graph is called \emph{self-bounding} if in each step: (i) it proceeds from the current vertex $i$ to a non-forbidden vertex $j \in \Gamma_i$, and (ii) adds to the set of forbidden vertices all neighbors of $i$ greater or equal to $j$.  We denote by $\mathcal{B}(G)$ the set of self-bounding walks on $G$.
\end{definition}

Self-bounding walks have been related before to the independent set polynomial \cite{SS, dror_IS}. Notably, in the seminal work of Weitz~\cite{dror_IS}, the vertices of the computation tree for $Z_G$ are the self-bounding walks on $G$ and the parent-child relationship amounts to the child extending the parent by one vertex. \medskip

Given a set of walks $\mathcal{W}$ and a walk $w \in \mathcal{W}$, we write $\Ext(w) 
$ for the set of walks in $\mathcal{W}$ that extend $w$ by one vertex. We let $\ter(w)$ denote the last, i.e., terminal, vertex of a walk $w$.

\begin{definition}
$\bvec{p} \in [0,1)^n$ is \emph{\good} for set of walks $\mathcal{W}$, if there exists 
$\RR:\mathcal{W} \mapsto [0,1)$ such that for all $w \in \mathcal{W}$, 
\begin{equation}\label{eq:WalkRatio}
p_{\ter(w)} \le \RR(w) \displaystyle{\prod_{y \in \Ext(w)}\left(1-\RR\left({y}\right)\right)} \enspace .
\end{equation}
\end{definition}


Recall that $\mathcal{S}(G) = \{\bvec{p} \in [0,1)^n: \avoid(\bvec{p}, G) > 0\}$.

 
\begin{restatable}{thm}{AbstractWalkLLLa}\label{thm:AbstractWalkLLLa}
$\bvec{p} \in \mathcal{S}(G) \Longleftrightarrow \bvec{p}$ is \good\ for $\mathcal{B}(G)$, the set of self-bounding walks of $G$.
\end{restatable} 

Observe that even though the form of~\eqref{eq:WalkRatio} is identical to the form of the LLL condition, i.e., to the form of a \emph{sufficient} condition for $\bvec{p} \in \mathcal{S}(G)$, Theorem~\ref{thm:AbstractWalkLLLa} is a \emph{characterization} of membership in $\mathcal{S}(G)$. In other words, inequality~\eqref{eq:WalkRatio} leaves nothing on the table: all relaxations of Theorem~\ref{thm:AbstractWalkLLLa} to more tractable sufficient conditions, i.e., to local lemmata such as the LLL, arise merely by ``fattening" the set $\mathcal{B}(G)$ to  ``nicer" sets, i.e., to sets for which establishing $\bvec{p}$ as \good\ can be done with a certificate whose size is only polynomial in $n$ (whereas, typically, $|\mathcal{B}(G)|$ is exponential in $n$). An essential element of such a collapse is to  pass to {\bf infinite} supersets of $\mathcal{B}(G)$, which we will then partition into a polynomial in $n$ number of equivalence classes. For example, as we will see shortly, the asymmetric LLL corresponds  to the set of \emph{all} walks on $G$.
\begin{definition}
Let $\sim$ be the {relation} where $w \sim w^\prime$ if $\ter(w) = \ter(w^\prime)$, and $\{z: wz \in \mathcal{W}\} = \{z: w^\prime z \in \mathcal{W}\}$.
We denote the equivalence class of a walk $w$ by $\widetilde{w}$ and the set of equivalence classes of a set of walks $\mathcal{W}$ by $C(\mathcal{W})$. Since functions $\ter$ and $\Ext$ are constant within each class, we allow ourselves to write $\ter(\widetilde{w})$ and $\Ext(\widetilde{w})$.
\end{definition}

\begin{restatable}{thm}{WalkRatioInequalityIff}\label{thm:WalkRatioInequalityIff}
For any set of walks $\mathcal{W}$ on $G$ and any $\bvec{p}\in [0,1)^n$, the following are equivalent:
\begin{enumerate}
\item
$\bvec{p}$ is \good\ for $\mathcal{W}$.
\item For every $\zz\in C(\mathcal{W})$ there exists $r_{\zz}\in [0,1)$ such that
$p_{\ter({\zz})} \le r_{\zz} \prod_{y \in \Ext(\zz)}(1-r_{\widetilde{y}})$.
\end{enumerate}
\end{restatable} 

Combining Theorems~\ref{thm:AbstractWalkLLLa} and~\ref{thm:WalkRatioInequalityIff} yields the following.
\begin{corollary}\label{cor:nyxta}
$\bvec{p} \in \mathcal{S}(G) \Longleftrightarrow$ There exists a set of walks $\mathcal{W} \supseteq \mathcal{B}(G)$ and $r_{\zz}\in [0,1)$ for each $\zz\in C(\mathcal{W})$, such that $p_{\ter({\widetilde{w}})} \le r_{\zz} \prod_{y \in \Ext({\widetilde{w})}}(1-r_{\widetilde{y}})$.
\end{corollary}

For example, the asymmetric LLL follows trivially from Corollary~\ref{cor:nyxta}, by taking $\mathcal{W}$ to be the set of \emph{all} walks on $G$. To see this, observe that when $\mathcal{W}$ is the set of all walks, $w \sim w'$ iff $\ter(w) = \ter(w')$, and thus the equivalence classes of $\Ext(w)$ correspond to the neighbors of $\ter(w)$. Thus, the condition ``for every $\zz\in C(\mathcal{W})$ there exists $r_{\zz}\in [0,1)$ such that $p_{\ter({\widetilde{w}})} \le r_{\zz} \prod_{y \in \Ext({\widetilde{w})}}(1-r_{\widetilde{y}})$," becomes ``for every $i \in [n]$, there exists $r_i \in [0,1)$ such that $p_i \le r_i \prod_{j \in \Gamma_i}(1-r_j)$," i.e., the asymmetric LLL condition.\smallskip

In fact, by taking $G = K_n$, Corollary~\ref{cor:nyxta} even recovers the \emph{union bound}. Specifically, it is easy to see that $\mathcal{B}(K_n)$ is the set of strictly descending sequences on $[n]$ and, by our Theorem~\ref{thm:AbstractWalkLLLa}, $\avoid(\bvec{p}, G) > 0$ iff $\bvec{p}$ is \good\ for the set of such sequences. A little bit of work then shows that $\bvec{p} \in [0,1)^n$ is \good\ for the set of descending sequences on $[n]$ iff $||\bvec{p}||_1 < 1$.

\subsection{A Hierarchy of Local Lemmata}

In order to keep the number of equivalence classes in Theorem~\ref{thm:WalkRatioInequalityIff} polynomial in $n$, the key is to consider walks for which the number of forbidden vertices at any moment is bounded by a quantity independent of $n$. For example, when we took $\mathcal{W}$ to be the set of all walks on $G$, this was trivially true since the set of forbidden vertices is always empty and, thus, there are as many equivalence classes as possible terminal vertices of a walk. Arguably, the next simplest class of walks containing $\mathcal{B}(G)$ are \emph{non-backtracking} walks, wherein the only forbidden vertex at any moment is the next-to-last vertex of the walk. Applying Corollary~\ref{cor:nyxta} when $\mathcal{W} \supseteq \mathcal{B}(G)$ is the set of non-backtracking walks, yields Theorem~\ref{thm:nbLLL}, below, which was first established in~\cite{ICALPIS} via an ad hoc inductive proof. Since for non-backtracking walks each equivalence class is represented either by a single vertex (corresponding to walks with a single vertex), or by the last and the next-to-last vertex of the walks (which, by definition, must be adjacent), we see that there is an equivalence class, and, thus, a parameter, for every vertex and every \emph{oriented edge} of $G$.

\begin{restatable}{thm}{nbLLL}\label{thm:nbLLL}\emph{(Non-backtracking LLL~\cite{ICALPIS})} 
Let $A = \{(i,j) : \{i,j\} \in E\}$. For $a=(i,j) \in A$, let $\Gamma_{a} = \{(j,k) \in A: k \neq i\}$. If there exist 
$0 \le \{r_i\}_{i \in [n]}, \{r_a\}_{a \in A}<1$ 
such that for all $i \in [n]$ and $a=(i,j) \in A$,
\[
p_i \le r_i  \prod_{j \in \Gamma_i} \left(1-r_{(j,i)}\right)\qquad \text{and} \qquad p_j \le r_a \prod_{y \in \Gamma_a}(1-r_y) \enspace ,
\]
then $\bvec{p} \in \mathcal{S}(G)$.
\end{restatable}
As shown in~\cite{ICALPIS}, Theorem~\ref{thm:nbLLL} implies than in the asymmetric LLL condition the expression $\prod_{j \in \Gamma_i} (1-r_j)$ can be replaced by $\prod_{j \in \Gamma_i} (1-r_j)/(1-r_i r_j)$, yielding a strict improvement of the asymmetric LLL. \medskip

Naturally, one can similarly apply Corollary~\ref{cor:nyxta} with $\mathcal{W}$ being the set of 2-non-backtracking walks, i.e., the walks where the set of forbidden vertices at any moment comprises the penultimate vertex and its predecessor, the set of 3-non-backtracking-walk, etc., and get a hierarchy of increasingly sharper local lemmata similar to Theorem~\ref{thm:nbLLL}. As it turns out, though, this is not ideal. This is because the set of forbidden vertices in self-bounding walks is, generally, significantly bigger than the set of  visited vertices, and giving up on those additional vertices is costly. Instead, we get much better results if we adopt the following viewpoint. We will define a collection of sets $S_1, \ldots, S_q \subseteq [n]$ and think of each set $S_i$ as a ``filter," in the following sense: $S_i$ will ``accept" walk $w$ only if every time $w$ enters $S_i$, the subwalk it performs within $S_i$ until it exits, belongs in $\mathcal{B}(G[S_i])$. We will only consider those walks accepted by \emph{all} filters. Formally, we have the following.
\begin{definition}\label{def:Csb}
For a walk $w$, every contiguous substring of $w$ is called a \emph{subwalk} of $w$. Given $S \subseteq [n]$, the $S$-\emph{restriction} of a walk $w$ is the set of subwalks of $w$ whose vertices are all in $S$. We say that $w$ is \emph{self-bounding within $S$} if every walk in the $S$-restriction of $w$ is in $\mathcal{B}(G[S])$. Given a family $\mathbf{S} = \{S_1, \ldots, S_q\}$ of subsets of $[n]$, we say that a walk is \emph{$\mathbf{S}$-self-bounding} if it is self-bounding within $S_i$ for every $i \in [q]$.
\end{definition}

Obviously, the main point here is that we are completely free to adapt our choice of $\mathbf{S}$ to the graph $G$, so that for the chosen complexity of the condition we control the fattening of $\mathcal{B}(G)$ as well as possible. At the same time, key to being able to certify that a vector $\bvec{p}$ is \good\ for a set of walks, will be that the set is infinite, capturing that the walks have ``limited memory." Naturally, it is not a priori clear that the filtering approach outlined above, which clearly allows for infinite walks, will lead to a finite (let alone polynomial) number of equivalence classes and, thus, a workable condition. Moreover, it is also not a priori obvious that the resulting set of walks will include all self-bounding walks, a necessary condition for the set to be useful. Theorem~\ref{thm:finiteclasses}, below, puts both of these concerns to rest.

\begin{restatable}{thm}{finiteclasses}\label{thm:finiteclasses}
Let  $\mathbf{S}$ be any family of subsets of $[n]$ and let $\mathcal{K}(\mathbf{S})$ be the set of \emph{$\mathbf{S}$-self-bounding} walks on $G$. $\mathcal{K}(\mathbf{S}) \supseteq \mathcal{B}(G)$ and $C(\mathcal{K}(\mathbf{S}))$ is finite.
\end{restatable}

\begin{example*}

\mbox{}

\begin{enumerate}[label=(\roman*)]
\item\label{LLL:asymm}
If $\mathbf{S} = \emptyset$, then $\mathcal{W}$ is the set of all walks on $G$ and we get the asymmetric LLL.
\item
If $\mathbf{S} = \{ \{i,j\} : \{i,j\} \in E \}$, then $\mathcal{W}$ is the set of non-backtracking walks on $G$ and we get Theorem~\ref{thm:nbLLL}.
\item\label{LLL:neigh}
If $\mathbf{S} = \{ \Gamma_i \cup \{i\} : i \in [n] \}$, we get a local lemma that dominates the \emph{cluster expansion} LLL~\cite{cluster}.
\item\label{LLL:shear}
If $\mathbf{S} = \{[n]\}$, then  $\mathcal{W} = \mathcal{B}(G)$ and we recover Shearer's criterion.
\end{enumerate}
\end{example*}

Between the local lemmata corresponding to cases~\ref{LLL:asymm}--\ref{LLL:neigh} above and the global condition of Shearer in case~\ref{LLL:shear}, one can consider ``intermediate" regimes with corresponding complexity. Prior to our work, the only effort in that direction was the \emph{decomposition local lemma} of Kolipaka, Szegedy, and Xu~\cite{cliqueLLL}, of which the Clique LLL (CLLL) is the more frequently applied special case. While~\cite{cliqueLLL} does not draw any connection to walks on the graph $G$, such a connection can be established with some effort and it illuminates why, for any given family of sets $\mathbf{S}$ (and thus two conditions of similar complexity), our method dominates  the decomposition LLL of~\cite{cliqueLLL}, often dramatically. Adopting the filter view discussed above, a walk on $G$ is ``accepted" in the method of~\cite{cliqueLLL}, i.e., contributes to the fattening, if it can be partitioned into subwalks each of which lies entirely within some filter $S_i$ and is self-bounding within $S_i$. While the two approaches are \emph{identical} when each edge of $G$ belongs in exactly one filter, as soon as any edge is included in multiple filters, our method acquires an advantage, as it is more selective in the walks it accepts. Since in typical applications, the filters correspond, roughly, to depth-$d$ vertex neighborhoods, the amount of filter overlap (and the corresponding gap in performance) blows up rapidly in $d$. Arguably, though, the bigger difference between the two approaches is that while our method enjoys monotonicity, i.e., enlarging filters and/or adding new ones, can only shrink the set of accepted walks and, thus, improve things, that is \emph{not} the case for the decomposition LLL of~\cite{cliqueLLL}. As a result, besides performing worse for any given set of filters, designing a good decomposition (set of filters) for the decomposition LLL is far from obvious. See, for example, the discussion in Section~\ref{sec:numbers} regarding the choice of decomposition for the cubic and hexagonal lattices.

\subsection{Submodular Functions and the Independent Set Polynomial}

Let $f : \{0,1\}^n \mapsto \Real_{\ge 0}$. We will find it useful to also think of $f$ as a function on the subsets of $[n]$, by thinking of each $\bvec{x} \in \{0,1\}^n$ as the characteristic vector of a subset of $[n]$. For $i \in [n]$, and $S \subseteq [n]\setminus \{i\}$, let 
\[
{\Delta_i{f(S)}:= f(S\cup \{i\}) - f(S)}, \text{ and } f(i\mid S) := \frac{f(S\cup \{i\})}{f(S)} \enspace .
\]
${\Delta_i{f(S)}}$ is the \emph{discrete derivative} of $f$ with respect to $i$ at $S$, measuring the \emph{additive} change of $f$ when moving along the direction of $i$, a fundamental quantity in the analysis of functions on the Boolean cube. $f(i\mid S)$ also measures the change of $f$ along the direction of $i$, but in a \emph{multiplicative} sense.
We call $f$ \emph{increasing/decreasing} if $\Delta_i f(S) \ge 0$ is positive/negative, respectively, for all $i, S$. We call $f$ \emph{supermodular}/\emph{submodular} if $\Delta_i{f}$ is increasing/decreasing, respectively, for all $i \in [n]$. Analogously, we call $f$ \emph{log-supermodular}/\emph{log-submodular}, if $f(i\mid S)$ is increasing/decreasing with respect to $S$, respectively, for all $i,S$. \medskip

We prove \emph{lower bounds for supermodular functions,} i.e., for the functions that capture the notion of ``diminishing returns" for decreasing functions  (since for decreasing functions $\Delta_i f$ is negative, in order for $|\Delta_i f|$ to be diminishing,  $\Delta_i f$ must be increasing, i.e., $f$ must be supermodular). 
To bound any function from below we need (a) an initial condition and (b) some control over its rate of decrease. In our setting, this will take the form (a) $f(\emptyset) = f(\bvec{0}) > 0$ and (b) lower bounds for \emph{some} of the marginal ratios $f(i|S)$. Analogously to the LLL, the information in (b) will be represented by a vector $\bvec{p} \in [0,1)^n$, and a graph $G$ on $[n]$. Specifically, for each $i \in [n]$, the neighbors of $i$ in $G$ will represent the dimensions for which we have no control whatsoever of the interaction with dimension $i$, whereas for any set of dimensions  $S \subseteq [n] - (\Gamma_i \cup \{i\})$  comprised of non-neighbors of $i$, the decrease of $f$ when $i$ is added to $S$ will be bounded \emph{multiplicatively} as
\begin{equation}\label{eq:promise}
f(S\cup \{i\}) \ge (1-p_i) f(S) \enspace. 
\end{equation}
\begin{definition}
A  function $f$ \emph{factorizes according to $\bvec{p}, G$} if~\eqref{eq:promise} holds for all $i \in [n]$ and $S \subseteq [n] - (\Gamma_i \cup \{i\})$.
\end{definition}

It is easy to see that the probabilistic LLL setting is a special case of the above.
\begin{example}\label{ex:only}
Let $\Omega$ be an arbitrary set with probability measure $\mu$, let $F_1, \ldots, F_n \subseteq \Omega$, and write $\mu(F_i) = p_i$. It is easy to see that $G$ is a dependency graph for the events corresponding to $F_1, \ldots, F_n$ iff $\mu$ factorizes according to $\bvec{p},G$. To see that $\mu$ is supermodular, observe that for every $i \in [n]$, and 
$S \subseteq T \subseteq [n] \setminus \{i\}$,
\[
\Delta_i\mu(S) = 
{\mu\left(\overline{F_i}\cap\bigcap_{j \in S}\overline{F_j}\right)}
-{\mu\left(\bigcap_{j \in S}\overline{F_j}\right)}
=
-{\mu\left({F_i}\cap\bigcap_{j \in S}\overline{F_j}\right)}
\le
-{\mu\left({F_i}\cap\bigcap_{j \in T}\overline{F_j}\right)}
=
\Delta_i \mu(T) \enspace.
\]
\end{example}\bigskip

We prove that for supermodular functions that factorize according to $\bvec{p},G$, the independent set polynomial of $G$ at $-\bvec{p}$ characterizes their maximum drop relative to $f(\emptyset)$.  For example, applying our Theorem~\ref{thm:ourShearer} to the function $\mu$ of Example~\ref{ex:only}, recovers exactly Shearer's criterion for the supermodular function $\avoid(\bvec{p}, G)$.
\begin{restatable}{thm}{ourShearer}\label{thm:ourShearer}
Let $G$ be a graph on $[n]$, and $\bvec{p}\in [0,1)^n$.
\begin{itemize}
\item[(i)] If $Z_G(-\bvec{p}; S)>0$, for every $S \subseteq [n]$ and $f$ is a supermodular function with $f(\bvec{0})>0$ that factorizes according to $G,\bvec{p}$, then for every $\bvec{x} \in \{0,1\}^n$,
\[
f(\bvec{x}) \ge f(\bvec{0}) \cdot Z_G(- p_1 x_1, -p_2 x_2, \ldots, - p_n x_n) > 0 \enspace .
\]
\item[(ii)] If $Z_G(-\bvec{p}; S) < 0$ for some $S \subseteq [n]$, then there exists a supermodular function $f$ with $f(\bvec{0})>0$ that factorizes according to $G,\bvec{p}$, such that $f(\bvec{x}) = 0$ for some $\bvec{x} \in \{0,1\}^n$.
\end{itemize}
\end{restatable}


\begin{corollary}[Supermodular Local Lemma]
If $f$ is a supermodular function with $f(\emptyset) > 0$ that factorizes according to $\bvec{p}, G$, and there exist $r_1, \ldots, r_n \in [0,1)$ such that $p_i \le r_i \prod_{j \in \Gamma_i} (1-r_j)$, for all $i \in [n]$, then 
\[
f(S) \ge f(\emptyset) \prod_{i\in S} (1-r_i) > 0 \enspace ,
\]
for every $S \subseteq[n]$.
\end{corollary}
\begin{proof}
Under the hypothesis, Dobrushin's condition implies $Z_G(-\bvec{p}; S) \ge \prod_{i\in S} (1-r_i) > 0$, for every $S \subseteq[n]$. Thus, Theorem~\ref{thm:ourShearer} implies $f(S) \ge f(\emptyset) \cdot Z_G(-\bvec{p}; S)$.
\end{proof}

%

\comment{What makes Theorem~\ref{thm:ourShearer} a significant generalization of Shearer's criterion is that, unlike the probabilistic setting, we do not even need to support the notion of the ``complement" for each of our $n$ sets, let alone require strong properties such as \emph{conditional complementarity}, i.e., $f(A \cap B)/f(B) +f(\overline{A} \cap B)/f(B)=1$. This great freedom yields quantum versions of all LLL results automatically, as we demonstrate next.} 

\subsection{Quantum Local Lemmata}
Given a vector space $V$ and a subspace $X \subseteq V$, let $R(X) := \frac{\dim (X)}{\dim(V)}$ be the relative dimension of $X$ with respect to $V$ and write $R(A \mid B) := \frac{R(A \cap B)}{R(B)}$. For $U \subseteq V$, let $X^{\perp U}$ denote the orthogonal complement of $X$ with respect to $U$, i.e., $\{v \in U: \langle v, x \rangle= 0, \text{for all } x \in X\}$.
\begin{definition}
A subspace $X$ is \emph{mutually $R$-independent} of subspaces $X_1, \ldots, X_{\ell}$ if $R(X \cap X_1 \cap \cdots \cap X_{\ell}) = R(X)\, R(X_1 \cap \cdots \cap X_{\ell})$. A graph $G$ on $[n]$ is an \emph{$R$-dependency graph} for subspaces $X_1, \ldots, X_n$ if for every $i \in [n]$ and every $S \subseteq [n] - (\Gamma_i \cup \{i\})$, the subspace $X_i$ is mutually independent of subspaces $\{X_j\}_{j\in S}$.
\end{definition}

A key difference between probabilistic independence and $R$-independence is that while for any two sets $A,B$, we have $\Pr(A \mid B) + \Pr(\overline{A} \mid B) =1$, it is \emph{not} generally the case that $R(A \mid B) + R(A^{\perp} \mid B) = 1$. Overcoming this problem was the main contribution of Ambainis et al.~\cite{QLLL}, who proved an analogue of the LLL, and of Sattath et al.~\cite{SPNAS} who proved an analogue of Shearer's criterion.\medskip

Let $X_1, \ldots, X_n$ be subspaces with $R$-dependency graph $G$ and $R(X_i) = 1 - p_i$, for all $i \in [n]$. \begin{restatable}{thm}{QLLL}\emph{(Quantum  LLL~\cite{QLLL})}\label{thm:QLLL}
If there exist  $r_1, \ldots, r_n \in [0,1)$ such that $p_i \le r_i \prod_{j \in \Gamma_i} (1-r_j)$ for all $i \in [n]$, then $R(\bigcap_{i=1}^n X_i) \ge \prod_{i\in[n]} (1-r_i) > 0$. 
\end{restatable}

\begin{restatable}{thm}{QuantumShearer}\emph{(Quantum Shearer~\cite{SPNAS})}\label{thm:QuantumShearer}
$\bvec{p} \in \mathcal{S}(G) \Longleftrightarrow R(\bigcap_{i=1}^n X_i) >0$.
\end{restatable}

Below we show how Theorems~\ref{thm:QLLL} and~\ref{thm:QuantumShearer} follow readily from Theorem~\ref{thm:ourShearer}.

\begin{proof}[Proof of Theorems~\ref{thm:QLLL} and~\ref{thm:QuantumShearer}]
For $S \subseteq [n]$, let $f(S) = R\left(\bigcap_{j \in S}X_j \right)$, where $f(\emptyset)  = 1$. We will prove that $f$ is supermodular and factorizes according to $\bvec{p},G$. The latter is trivial as for all $i \in [n]$ and $S \subseteq [n] - (\Gamma_i \cup \{i\})$,
\[
f(i \mid S) = \frac{f(S \cup \{i\})}{f(S)}
 = \frac{R\left({X_i}\cap\bigcap_{j \in S}{X_j}\right)}
{R\left(\bigcap_{j \in S}{X_j}\right)}
 = \frac{R(X_i)R(\bigcap_{j\in S} X_j)}
{R\left(\bigcap_{j \in S}{X_j}\right)}
= R(X_i)
= 1- p_i \enspace.
\]
Regarding supermodularity we observe that for every $i \in [n]$ and $S \subseteq [n] - \{i\}$,
$$
\Delta_if(S) 
= R(S \cup \{i\}) - R(S) 
	= - \frac{\dim\left(X_i^{\perp (\bigcap_{j \in S} X_j)}\right)}{\dim(V)} \enspace .
$$
Therefore, for every $i \in [n]$ and every $S \subseteq T \subseteq [n]\setminus \{i\}$, we have $\Delta_if(S) \le \Delta_if(T)$. \smallskip

Since $f(\emptyset) =1 >0$, $f$ is supermodular, and $f$ factorizes according to $\bvec{p},G$, we can apply Theorem~\ref{thm:ourShearer} and get $f([n]) = R(\bigcap_{i=1}^n X_i) >0$ iff $Z_G(-\bvec{p};S)>0$ for every $S \subseteq [n]$, i.e., Theorem~\ref{thm:QuantumShearer}. Moreover, if  $\bvec{p} \in \mathcal{S}(G)$, then $f([n]) = R(\bigcap_{i=1}^n X_i) = f(\bvec{1}) \ge f(\bvec{0}) Z_G(-\bvec{p}) = Z_G(-\bvec{p})$. Theorem~\ref{thm:QLLL} now follows by observing that under its hypothesis, Dobrushin's condition implies that $\bvec{p} \in \mathcal{S}(G)$ and, in particular, $Z_G(-\bvec{p})\ge \prod_{i\in[n]} (1-r_i)$.
\end{proof}

\section{Proof of Theorem~\ref{thm:AbstractWalkLLLa}}

It will be helpful to introduce the quantity $\ratio(\bvec{p}; (i, S)) := 1-Z(i \mid S) =\ratio(i, S) $ and rewrite~\eqref{eq:identical} as
\begin{equation}\label{eq:ratioRec}
\ratio(i, S) = p_i \prod_{\ell=1}^{d} \frac{1}{1- \ratio(j_\ell, S_\ell)} \enspace .
\end{equation}

\begin{lemma}\label{lem:StoR}
$\bvec{p} \in \mathcal{S}(G) \Longleftrightarrow$  $\ratio(\bvec{p}; (i,S)) < 1$, for every $i \in [n]$, and $S \subseteq [n]\setminus\{i\}$.
\end{lemma}

\begin{proof}

For arbitrary $S = \{j_1, \ldots, j_q\} \subseteq [n]$, we write $S_{\ell} = S - \{j_1, \ldots, j_{\ell}\}$, with $S_0 =S$. By telescoping, $Z(S) = \prod_{\ell \in [q]} Z(S_{\ell-1})/Z(S_{\ell}) = \prod_{\ell \in [q]} Z(j_{\ell} \mid S_{\ell})= \prod_{\ell \in [q]} (1 - \ratio(j_\ell, S_\ell))$.\smallskip

\noindent $\Longleftarrow$ 
If $\ratio(\bvec{p}; (i,S))<1$, for every $i \in [n]$, and $S \subseteq [n]\setminus\{i\}$, then the product $\prod_{\ell \in [q]} (1 - \ratio(j_\ell, S_\ell)) = Z(S)$ is strictly positive for every $S \subseteq [n]$ which, by Shearer's criterion, implies $\bvec{p} \in \mathcal{S}(G)$.\smallskip

\noindent $\Longrightarrow$ 
We will prove the contrapositive, i.e., that if there is $i \in [n]$ and $S \subseteq [n] - \{i\}$ such that $\ratio(\bvec{p}; (i,S)) \ge 1$, then $\bvec{p} \not\in \mathcal{S}(G)$. Let $\lambda = \min\{\theta >0 : \text{there exist $i \in [n]$ and $S \subseteq [n] - \{i\}$ such that $\ratio(\theta\bvec{p}; (i,S)) \ge 1$}\}$. {By continuity,} there exist  $i \in [n]$ and $S \subseteq [n] - \{i\}$ such that $\ratio(\lambda \bvec{p}; (i,S)) = 1$. In that case, though, $Z(\lambda \bvec{p}; S \cup \{i\}) = 0$, which by Shearer's criterion implies $\bvec{p} \not\in \mathcal{S}(G)$.
\end{proof}

\begin{lemma}\label{eq:SubWalkMon}
If $\bvec{p}$ is \good\ for $\mathcal{W}$, then $\bvec{p}$ is also \good\ for every $\mathcal{W'} \subseteq \mathcal{W}$.
\end{lemma}
\begin{proof}
Since $\bvec{p}$ is \good\ for $\mathcal{W}$, there exists $L : \mathcal{W} \mapsto [0,1)$ satisfying~\eqref{eq:WalkRatio}. Since $L(w) < 1$ for all $w \in \mathcal{W}$, removing any $w$ from $\mathcal{W}$ can only increase the r.h.s.\ of~\eqref{eq:WalkRatio} for the remaining walks.
\end{proof}

\begin{definition}\label{def:lhat}
For finite $\mathcal{W}$, let $\widehat{L} := \widehat{L}_{\bvec{p}, \mathcal{W}} : \mathcal{W} \mapsto [0,1)$ be the function defined by the recurrence 
\begin{equation}\label{eq:madhat}
\widehat{\RR}(w)  =  p_{\ter(w)} \prod_{y \in \Ext(w)} \frac{1}{1-\widehat{\RR}(y)} \enspace .
\end{equation}
\end{definition}

\begin{lemma}\label{eq:IntoEq}
If $L : \mathcal{W} \mapsto [0,1)$ satisfies~\eqref{eq:WalkRatio} for all $w \in \mathcal{W}$, then $L(w) \ge \widehat{L}(w)$, for all $w \in \mathcal{W}$.
\end{lemma}

\begin{proof}
Let $L:\mathcal{W} \mapsto [0,1)$ be an arbitrary function satisfying \eqref{eq:WalkRatio}. Since $\mathcal{W}$ is finite,  $\widehat{L}$ is well defined as a formal expression of $\bvec{p}$, and thus, all it remains to show is that $\widehat{L}(w) \le {L}(w)$, for all $w \in \mathcal{W}$. Let $\{w_1, w_2, \ldots, w_\ell\}$ be an ordering of $\mathcal{W}$ such that  for all $i \in [\ell]$, $w_i$ has no extensions in $\{w_i, \ldots, w_\ell\}$ (such an ordering must exist since $\mathcal{W}$ is finite). We show that $\widehat{L}(w_i) \le L(w_i)$, using induction on $i$. Let ${w_1, \ldots, w_s}$ be the walks in $\mathcal{W}$ having no extensions, then trivially, $\widehat{L}(w_j)= p_{\ter(w_j)} \le L(w_j)$, for all $j \in [s]$. Assume now we have established $\widehat{L}(w_j) \le L(w_j)$, for all $j<i$, and let $\Ext(w_i) = \{j_1, \ldots, j_d \} \subseteq [i-1]$. Then 
\[
\widehat{L}(w_i)
=\frac{p_{\ter(w_i)}  }{\displaystyle\prod_{k \in [d]}\left(1-\widehat{\RR}\left({w_{j_k}}\right)\right)}
\le \frac{p_{\ter(w_i)}  }{\displaystyle\prod_{k \in [d]}\left(1-{\RR}\left({w_{j_k}}\right)\right)} \le {L}(w_i)\enspace,
\]
where the inequality follows from the inductive hypothesis applied on $w_{j_1}, \ldots, w_{j_d}$, and the fact that the function $x \mapsto 1/(1-x)$ is increasing in $[0,1)$, while the second inequality  follows from our hypothesis.
\end{proof}

We restate Theorem~\ref{thm:AbstractWalkLLLa} below for convenience.

\AbstractWalkLLLa*

\begin{proof}[Proof of Theorem~\ref{thm:AbstractWalkLLLa}] 
Recall that a walk is called self-bounding if in each step: (i) it proceeds from the current vertex $i$ to a non-forbidden vertex $j \in \Gamma_i$, and (ii) adds to the set of forbidden vertices, $\mathcal{F}$, all neighbors of $i$ greater or equal to $j$. Recall that we denote the set of self-bounding walks on a graph $G$ by $\mathcal{B}(G)$. For $w \in \mathcal{B}(G)$, let $\mathcal{F}(w)$ denote the set of forbidden vertices of $w$, i.e., the set $\mathcal{F}$ when we reach $\ter(w)$.\medskip

\noindent $\Longrightarrow$
Let $\mathcal{W} = \mathcal{B}(G)$. For $w \in \mathcal{W}$, let
\[
L(w) = \ratio	\left(
				\bvec{p} ; \left(	
								\ter(w), [n]\setminus{\mathcal {F}}(w)
						\right)
			\right) \enspace .
\]
By unfolding recurrence~\eqref{eq:ratioRec}, it is not hard to see that $L$ satisfies~\eqref{eq:WalkRatio} as an equality. Thus, we are left to prove that $L < 1$. This  follows from Lemma~\ref{lem:StoR}, as $\bvec{p} \in \mathcal{S}(G) \Longrightarrow \ratio(i,S) < 1$, for every $i \in [n]$ and $S \subseteq [n]\setminus\{i\}$. \smallskip

\noindent $\Longleftarrow$
We will prove that if $\bvec{p}$ is \good\ for $\mathcal{B}(G)$, then $\ratio(i,S) < 1$ for every $i \in [n]$ and $S \subseteq [n]\setminus\{i\}$, which, by Lemma~\ref{lem:StoR}, implies $\bvec{p} \in \mathcal{S}(G)$.  Let $i \in [n]$ and $S \subseteq [n]\setminus\{i\}$ be arbitrary and write $S^+ := S\cup\{i\}$.


Write $(i)$ for the trivial path consisting of vertex $i$ alone. It is easy to check that recurrence~\eqref{eq:ratioRec} for $\ratio\left(\bvec{p};(i ,S)\right)$ coincides with recurrence~\eqref{eq:WalkRatio} for $\widehat{\RR}_{\mathcal{B}(G[S^+])}(\bvec{p};(i))$ and, thus,
\[
\ratio\left(\bvec{p};(i ,S)\right) = \widehat{\RR}_{\mathcal{B}(G[S^+])}(\bvec{p};(i)) \enspace .
\]

Now, the fact that $\bvec{p}$ is \good\ for $\mathcal{B}(G)$ implies that there exists $L:\mathcal{B}(G) \mapsto [0,1)$ satisfying~\eqref{eq:WalkRatio}, for every $w \in \mathcal{B}(G)$. Therefore, for every $i \in [n]$ and $S \subseteq [n] - \{i\}$, we can conclude that
\[
\ratio\left(\bvec{p};(i ,S)\right)= \widehat{\RR}_{\mathcal{B}(G[S^+])}(\bvec{p};(i)) \le \widehat{\RR}_{\mathcal{B}(G)}(\bvec{p};(i)) < L(\bvec{p};(i)) <1
\]
where the first inequality is due to Lemma~\ref{eq:SubWalkMon}, the second due to Lemma~\ref{eq:IntoEq}, and the last due to $L < 1$.
\end{proof}

\section{Proof of Theorem~\ref{thm:WalkRatioInequalityIff}}

Let $|w|$ denote the number of vertices of a walk $w$.
 
\begin{definition}
 For a set of walks $\mathcal{W}$, $k \in \Natural$, and $\bvec{p} \in [0,1)^n$, let $L_k: \mathcal{W} \mapsto \Real$ be defined as
\begin{equation}\label{eq:defRell}
\RR_k(\bvec{p};w) = 
\begin{cases}
0 																		&, \enspace |w| > k \\
\\
\displaystyle p_{\ter(w)}\prod_{z \in \Ext(w)} \frac{1}{1-\RR_k\left(\bvec{p};{z}\right)} 	&, \enspace |w| \le k \enspace .
\end{cases}
\end{equation} 
\end{definition}


\begin{lemma}\label{lem:mits}
If $\RR_{k}(\bvec{p} ; w) \le \RR_{k+1}(\bvec{p} ; w) < 1$ for all $k \in \Natural$ and all $w \in \mathcal{W}$, then $\bvec{p}$ is \good\ for $\mathcal{W}$.
\end{lemma} 
\begin{proof}
Since $L_k(\bvec{p};w)$ is increasing in $k$ and bounded, it follows that $\lim_{k \to \infty} L_k(\bvec{p};w) =\widehat{L}(\bvec{p} ; w) = \widehat{L}(w)$, where $\widehat{L}$ is per Definition~\ref{def:lhat}. Moreover, $\widehat{L}$ satisfies~\eqref{eq:WalkRatio} as an equality, establishing that $\bvec{p}$ is \good\ for $\mathcal{W}$.
\end{proof}

We restate Theorem~\ref{thm:WalkRatioInequalityIff} below for ease of reference.

\WalkRatioInequalityIff*

\begin{proof}[Proof of Theorem~\ref{thm:WalkRatioInequalityIff}]

\mbox{}

\noindent
$(1 \implies 2$) The fact that $\bvec{p}$ is \good\ for $\mathcal{W}$ implies that there exists $L$ satisfying~\eqref{eq:WalkRatio} for all $w \in \mathcal{W}$. Observe that if $w \sim z$, their corresponding inequalities are identical other than, potentially, for the values $L(w), L(z)$. Setting $L(w) = L(z) = \min\{L(w), L(z)\}$ retains the validity of the two inequalities, while all other inequalities are (``at least as") valid as before. Therefore, if for each $\widetilde{w} \in C(\mathcal{W})$ we set $r_{\widetilde{w}}$ to the minimum value of $L$ over the class, the inequality $p_{\ter({\zz})} \le r_{\zz} \prod_{y \in \Ext(\zz)}(1-r_{\widetilde{y}})$ will hold. \medskip

\noindent
$(2 \implies 1$) 
Let $P(\ell, k)$ be the proposition: ``for all $w \in \mathcal{W}$ of length $\ell$ it holds that $L_{k-1}(w) \le L_{k}(w) \le r_{\tilde{w}}$". Since $r_{\zz} < 1$, by Lemma~\ref{lem:mits}, it suffices to prove $P(\ell, k)$, for all $\ell\ge 1$, and all $k \ge 2$. First, we observe that if $\ell > k$, then $L_{k-1}(w) = L_{k}(w) = 0 \le r_{\tilde{w}}$, for all $w$ with $|w| = \ell$. 

To show $P(\ell, k)$ in the regime $\ell \le k$, we use induction on $(k - \ell)$.  For $k=\ell$ we have that $L_{k-1}(w) = 0$, while $L_{k}(w) = p_{\ter(\tilde{w})} \le r_{\tilde{w}}$, for all $w$ with $|w| = \ell$. Assume now we have established $P(\ell ,k)$ for all $\ell, k$ with $k-\ell < s$, for some $s \ge0$. Let $w \in \mathcal{W}$ be an arbitrary word of length $\ell$, and let $k =s +\ell$. Then

\begin{equation}\label{eq:PrwtoIn}
L_{k-1}(w) 
= \frac{p_{\ter(w)}}{\displaystyle{\prod_{y \in \Ext(w)}\left(1-\RR_{k-1}(y)\right)}} 
\le \frac{p_{\ter(w)}}{\displaystyle{\prod_{y \in \Ext(w)}\left(1-\RR_{k}(y)\right)}}
= L_{k}(w) \enspace,
\end{equation}
and
 \begin{equation}\label{eq:DeuteroIn}
 L_{k}(w) =
\frac{p_{\ter(w)}}{\displaystyle{\prod_{y \in \Ext(w)}\left(1-\RR_{k}(y)\right)}} 
\le \frac{p_{\ter(w)}}{\displaystyle{\prod_{y \in \Ext(w)}\left(1-r_{\tilde{y}}\right)}}
\le r_{\tilde{w}}
\enspace,
\end{equation}
where the first inequalities in \eqref{eq:PrwtoIn}, \eqref{eq:DeuteroIn} follow from our inductive hypothesis $P(\ell+1,k)$, and the fact that the function $x \mapsto 1/(1-x)$ is increasing in $[0,1)$, while the second inequality of \eqref{eq:DeuteroIn} follows from our hypothesis.

\end{proof}

\section{Proof of Theorem~\ref{thm:finiteclasses}}
\comment{
We begin by observing that if $\mathbf{S}$ is any family of subsets of $[n]$ and $\mathbf{S}^\prime$ is obtained by  adding more subsets to $\mathbf{S}$, or by taking the union of some members of  $\mathbf{S}$, then every $\mathbf{S}^\prime$-self-bounding walk is $\mathbf{S}$-self-bounding.
Therefore, if $\mathbf{S} = \{S_1, \ldots, S_q\}$ and $\mathbf{S'}$ consists of the set $[n] = \overline{S_1} \cup S_1 \cup \cdots \cup S_q$, then $\mathcal{K}(\mathbf{S}) \supseteq \mathcal{K}(\mathbf{S'}) = \mathcal{K}(\{[n]\}) = \mathcal{B}(G)$.}
The fact that $\mathcal{K}(\mathbf{S})\supseteq \mathcal{B}(G)$, follows trivially from the observation that a walk violating self-boundness within $S$, for some $S\subseteq [n]$, also violates self-boundiness within $G$.

To show that $C(\mathcal{K})$ is finite, consider the following procedure for generating $\mathcal{K}(\mathbf{S})$. We start with $A = \emptyset$. First, we add to $A$ all walks of length one, i.e., all vertices of $G$. Then, for a walk $w \in A$ and $x \in [n]$, we add $wx$ to $A$ iff for every $S_i \in \mathbf{S}$, the greatest suffix of $wx$ that is fully contained in $S_i$ belongs in $\mathcal{B}(G[S_i])$. It is clear that when this process terminates $A =  \mathcal{K}(\mathbf{S})$ and that the total number of walks that can be added to $A$ in the course of the process is bounded by $\prod_{S_i \in \mathbf{S}}\left|\mathcal{B}(G[S_i])\right|$.

\section{Proof of Theorem~\ref{thm:ourShearer}}

We restate Theorem~\ref{thm:ourShearer} below for convenience.

\ourShearer*

\begin{proof}
\comment{\sout{For $S \subseteq [n]$, let $\bvec{p}_S$ denote the mutation of $\bvec{p}$ where all variables outside $S$ equal 0. Let $Z_G(\bvec{p};S)$ denote $Z_G(\bvec{p}_S)$.}} To simplify notation, below, we drop the dependence on $G,\bvec{p}$. For $i \in [n]$ and $S \subseteq [n]\setminus \{i\}$ define
\[
Z(i\mid S) := \frac{Z_G(-\bvec{p};  S\cup \{i\})}{Z_G(-\bvec{p}; S)} \enspace .
\]
\noindent(i)
We will prove that $f(i\mid S) 
\ge Z\left(i\mid S\right)$, for all $i \in [n]$ and $S\subseteq [n]\setminus \{i\}$. To see that this suffices, let $S \subseteq [n]$ be arbitrary, let $j_1, j_2, \ldots, j_q$ be any ordering of $S$, and for $\ell \in [q]$, write $S_{\ell} = S - \{j_1, \ldots, j_{\ell}\}$. Then, 
\[
\frac{f(S)}{f(\emptyset)} = \frac{f(S)}{f(S_1)} \frac{f(S_1)}{f(S_2)} 
\cdots \frac{f(S_{q-2})}{f(S_{q-1})} \frac{f(S_{q-1})}{f(\emptyset)} = \prod_{\ell \in [q]} f(j_{\ell} \mid S_{\ell}) \ge \prod_{\ell \in [q]} Z(j_{\ell} \mid S_{\ell}) = {Z_{G[S]}(-\bvec{p})}> 0 
 \enspace .
\]

To estimate $Z(i\mid S)$ we observe that the contribution to $Z(S \cup\{i\})$ of the sets including vertex $i$ equals $-p_i$ times the contribution of the sets not including $\Gamma^+_i$. Therefore, 
\begin{equation}\label{eq:papa}
Z(S\cup\{i\}) = Z(S) -p_{i}Z(S\setminus\Gamma_i) \enspace .
\end{equation}
With the above in mind, let $\{j_1, \ldots, j_d\}$ be an ordering of $\Gamma_i \cap S$, and write $S_\ell = S - \{j_1, \ldots, j_\ell\}$. Dividing
~\eqref{eq:papa} by $Z(S)$ and writing the ratio $Z(S\setminus\Gamma_i)/Z(S)$ in telescopic form yields
\begin{equation}\label{eq:identical}
Z(i\mid S) =  1 -p_{i} \frac{1}{\displaystyle{\frac{Z(S)}{Z(S\setminus\Gamma_i)}}}
=1 -p_{i} \frac{1}{\displaystyle{\prod_{\ell \in [d]}\frac{Z(S\setminus \{j_1, \ldots, j_{\ell-1}\} )}{Z(S\setminus \{j_1, \ldots, j_\ell\})}}}= 1 -p_{i} \prod_{\ell \in [d]}\frac{1}{\displaystyle Z(j_\ell \mid S_\ell)} \enspace .
\end{equation}

We are now ready to prove that $f(i\mid S) \ge Z\left(i\mid S\right)$, for all $i \in [n]$ and $S\subseteq [n]\setminus \{i\}$ by induction on $|S|$. Assume the claim holds for every $i \in [n]$ and every proper subset of $S$. Let $j_1, \ldots, j_d$ be an ordering of $S \cap \Gamma_i$, and write $S_\ell = S - \{j_1, \ldots, j_\ell\}$. Writing $f(i\mid S) = 1 + \frac{\Delta_if(S) }{f(S)}$, the first inequality in~\eqref{eq:split1} follows from the supermodularity of $f$, while 
the second follows from the fact that $f$ factorizes according to $G, \bvec{p}$,
\begin{equation}\label{eq:split1}
f(i \mid S) = 1 + \frac{\Delta_if(S) }{f(S)} 
\ge 1 + \frac{\Delta_if(S\setminus \Gamma_i)}{f(S)} 
= 1 + \left[f\left(i \mid S \setminus \Gamma_i\right) - 1\right] \frac{f(S \setminus \Gamma_i)}{f(S)}  
\ge 1 - p_i\frac{f(S \setminus \Gamma_i)}{f(S)} \enspace .
 \end{equation}
Next, we write $\frac{f(S)}{f(S\setminus \Gamma_i)} =
\frac{f(S)}{f(S_1)} \frac{f(S_1)}{f(S_2)} \cdots  \frac{f(S_{d-2})}{f(S_{d-1})} \frac{f(S_{d-1})}{f(S_d)} = 
\prod_{\ell \in [d]} f(j_\ell \mid S_\ell)$, since $S_d = S \setminus \Gamma_i$.  Thus, the r.h.s.\ of~\eqref{eq:split1} equals the l.h.s. of~\eqref{eq:split2}, the inequality in~\eqref{eq:split2} follows from the inductive hypothesis, while the equality from~\eqref{eq:identical}. 
\begin{equation}\label{eq:split2}
 1 - p_i\prod_{\ell=1}^{d}\frac{1}{f(j_\ell \mid S_\ell)} 
\ge 1 - p_i\prod_{\ell=1}^{d}\frac{1}{Z(j_\ell \mid S_\ell)} 
= Z(i\mid S)\enspace .
\end{equation}
\medskip

\noindent(ii)
\comment{
\red{To prove that (i) implies (ii) we prove the contrapositive}, i.e., that if $Z_G(-\bvec{p}; S) \le 0$ for some $S \subseteq [n]$, then we can find $f : \{0, 1\}^{n} \mapsto \Real_{\ge0}$ that is supermodular and factorizes according to $G, \bvec{p}$, such that $f([n]) = 0$. To do this, 
}
Let $\lambda = \min\{\theta >0 : \text{$\exists S \subseteq [n]$ such that $Z_{G[S]}(-\theta \bvec{p}) \le 0$}\}$ and let $f = Z_G(-\lambda \bvec{p})$. Thus, $f \ge 0$, since otherwise the continuity of the independent set polynomial would imply that $\lambda$ is not minimal. Readily $f$ inherits factorization according to $G, \bvec{p}$ from $Z_G$. To prove that $f = Z_G(-\lambda \bvec{p})$ is supermodular
we observe that, per~\eqref{eq:papa}, 
\begin{equation}\label{eq:lastone}
Z(-\lambda \bvec{p} ; S\cup\{i\}) - Z(-\lambda \bvec{p} ; S) = 
-\lambda p_{i} Z(-\lambda \bvec{p} ; S\setminus\Gamma_i) \enspace .
\end{equation}
Since the l.h.s.\ of~\eqref{eq:lastone} equals $\Delta_if(S)$ and $Z_G(-\lambda \bvec{p}) \ge 0$, we see that $f$ is decreasing.
Now, let $S \subseteq T$. Per~\eqref{eq:lastone}, we see that $\Delta_if(T) - \Delta_if(S) = 
-\lambda p_{i} Z(-\lambda \bvec{p} ; T\setminus\Gamma_i) + 
\lambda p_{i} Z(-\lambda \bvec{p} ; S\setminus\Gamma_i)$ which, since $f$ is decreasing, is positive. Thus, $\Delta_if(S)$ is increasing, i.e., $f$ is supermodular.
\comment{
Finally, \red{we observe that if (ii) holds}, then we can use telescoping to bound $f(\bvec{x})$ from \red{below} as
\[
f(\bvec{x}) = f(\bvec{0}) \prod_{i: x_i = 1} f(i \mid \{j < i : x_j = 1\}) \red{\le}  
 f(\bvec{0}) \prod_{i: x_i = 1} Z(i \mid \{j < i : x_j = 1\}) 
= f(\bvec{0}) Z_G(- p_1 x_1, \ldots, - p_n x_n) \enspace.
\]}
\end{proof}

\section{Negative-Fugacity Singularity for Hard-Core Model on Lattices}\label{sec:numbers}

When all variables of the independent set polynomial, i.e., all vertices of the graph, take the same value $\lambda$ (known in statistical physics as ``activity") the independent set polynomial becomes univariate, i.e., $Z_G(\lambda)$. Given a locally-finite, countable, pseudo-transitive graph $G$, typically a regular lattice such as $\mathbb{Z}^d$ our goal is to determine $\lambda_c = \sup\{\lambda: Z_G(-\theta \lambda) > 0, \text{ for all $\theta \in [0,1]$}$\} 
This number, $\lambda_c$, known as the  ``negative-fugacity singularity of the hard-core model," has been extensively studied for different lattices in the statistical physics community, primarily motivated by the Lee-Yang~\cite{leeyan52} approach to the study of phase transitions which implies that for $\lambda \in [0, \lambda_c)$ the hard-core model, i.e., the probability distribution on independent sets where each set $I \in \mathrm{Ind}(G)$  has probability proportional to $\lambda^{|I|}$, does not exhibit any phase transition. Very recently, Regts~\cite{regts2021absence} strenghtened the absence of phase transitions by showing that  for $\lambda \in [0, \lambda_c)$, the hard-core model exhibits decay of correlations in the form of \emph{strong spatial mixing}, a very helpful property for designing efficient algorithms for approximating the independent set polynomial for such $\lambda$.

Here we consider three important lattices: the square lattice ($\mathbb{Z}^2$), the cubic lattice ($\mathbb{Z}^3$), and the hexagonal planar lattice. Although the underlying graphs are infinite, due to their transitive structure, for any given integer $\ell$ we can find finite families of finite sets $\mathbf{S}=\mathbf{S}(\ell)$, such that the $\mathbf{S}$-self-bounding walks from the origin of length at most $\ell$, when extended periodically, generate a set that includes all self-bounding walks on the lattice. The lower bound we report below for each lattice, corresponds to choices for $\ell, \mathbf{S}$ for which the corresponding set of walks contains approximately one million walks. \smallskip

As can be seen from Table~\ref{tab:comp}, our method dominates all other known methods, for each lattice matching the conjectured value for $\lambda_c$ (derived using sophisticated but non-rigorous methods from statistical physics), to the third decimal digit. In particular, we note that that~\cite{cliqueLLL}, only reports a lower bound for $\lambda_c$ for the square lattice. To make the comparison more comprehensive we derived suitable decompositions ourselves, as follows. For the cubic lattice imagine the unit cubes of $\mathbb{Z}^3$ colored black and white in a 3D-chessboard pattern. Taking each set (filter) in the decomposition to be the set of vertices of a white unit cube, so that every edge is included in exactly one filter, gave $0.0695$. A more complex decomposition, yielded $0.0702$. For the hexagonal lattice, the number we report corresponds to the decomposition where each filter amounts to a single edge of the lattice, making the method equivalent to the non-backtracking LLL. We did this because several natural choices of larger filters only made things worse, manifesting the non-monotonicity of the decomposition LLL. For example, taking the filters to be the hexagonal faces gives only $0.067$, much less than taking the filters to be individual edges which gives $0.1481$.

\begin{table}[ht]
\centering
\begin{tabular}[t]{l|cccccc}
						& Asymmetric 				& Cluster 			& Non-back 				& Decomposition  	& This work 	& Numerical\\
						& \cite{erdHos1991lopsided}	& \cite{cluster}		& \cite{ICALPIS}			& \cite{cliqueLLL} 	& 			& \cite{todo1999transfer}\\
\hline
Square ($\mathbb{Z}^2$)
						& 0.0819 				& 0.0896 			& 0.1054					& 0.1130 		& 0.1191 	& 0.1193 \\
						\hline
Cubic ($\mathbb{Z}^3$)
						& 0.0566					& 0.0601				& 0.0669					& 0.0702		& 0.0743		& 0.0744\\
						\hline
Hexagonal
						& 0.1054					& 0.1190				& 0.1481					& 0.1481 		& 0.1542		& 0.1547\\
\hline
\end{tabular}
\caption{Lower bounds for $\lambda_c$ for different lattices and methods.}\label{tab:comp}
\end{table}

Finally, we note that the astute reader may have observed that the lower bound for the square lattice given by the non-backtracking LLL coincides with the lower bound given for the \emph{hexagonal} lattice by the asymmetric LLL. Remarkably, this is \emph{not} a coincidence. In the eyes of the non-backtracking LLL, each vertex of the square lattice (other than the origin) has degree three, since the predecessor of a vertex along a walk cannot be reused. On the other hand, in the hexagonal lattice the asymmetric LLL can use all three neighbors of each vertex, as it corresponds to the set of all walks.

\newpage
\bibliographystyle{alpha}
\bibliography{Papers} 
\end{document}